\documentclass[12pt]{amsart}
\usepackage{amssymb}
\usepackage{amscd}
\usepackage[all]{xypic}

\newcounter{alphalistcntr}

\theoremstyle{plain}
\newtheorem{theorem}{Theorem}[section]
\newtheorem{lemma}[theorem]{Lemma}
\newtheorem{corollary}[theorem]{Corollary}
\newtheorem{prop}[theorem]{Proposition}

\theoremstyle{remark}

\newtheorem{remark}[theorem]{Remark}

\newtheorem*{note*}{Note}
\newtheorem*{remark*}{Remark}
\newtheorem*{example*}{Example}

\theoremstyle{definition}
\newtheorem*{definition*}{Definition}
\newtheorem{definition}[theorem]{Definition}

\newcommand{\Z}{\mathbb{Z}}

\newcommand{\Q}{\mathbb{Q}}
\newcommand{\C}{\mathbb{C}}
\newcommand{\N}{\mathbb{N}}

\newcommand{\Gal}{\mathrm{Gal}}
\newcommand{\tensor}{\otimes}

\newcommand{\disc}{\mathrm{Disc}}
\newcommand{\Norm}{\mathrm{Norm}}

\newcommand{\tr}{\mathrm{Tr}}
\renewcommand{\mod}{\mathrm{mod}}

\title[Relative Galois module structure of absolutely abelian number fields]{Relative Galois module structure of rings of integers of absolutely abelian number fields}
\author{Henri Johnston}
\thanks{The author was partially supported by the Jersey Scholarship
from the States of Jersey Education, Sport and Culture Committee, 
and by a grant from the Deutscher Akademischer Austausch Dienst.}
\address{Dept. of Math. \\ Malott Hall \\ Cornell University \\
Ithaca, NY 14853-4201 \\ USA}
\curraddr{Fakult\"at f\"ur Informatik\\
Institut f\"ur theoretische Informatik und Mathematik\\
Universit\"at der Bundeswehr M\"unchen\\
85577 Neubiberg\\
Germany}
\email{henri@math.cornell.edu}
\urladdr{http://www.math.cornell.edu/$\sim$henri}
\date{5th July 2007}

\begin{document}

\begin{abstract}
Let $L/K$ be an extension of number fields where $L/\Q$ is abelian. We define such
an extension to be Leopoldt if the ring of integers $\mathcal{O}_{L}$ of $L$ is free
over the associated order $\mathcal{A}_{L/K}$. Furthermore we define an abelian
number field $K$ to be Leopoldt if every finite extension $L/K$ with $L/\Q$ abelian
is Leopoldt in the sense above. Previous results of Leopoldt, Chan \& Lim, Bley, and Byott \& Lettl
culminate in the proof that the $n$-th cyclotomic field $\Q^{(n)}$ is Leopoldt for
every $n$. In this paper, we generalize this result by giving more examples
of Leopoldt extensions and fields, along with explicit generators.
\end{abstract}

\maketitle

\section{Introduction}\label{intro}

Let $L/K$ be a finite Galois extension of number fields with Galois group $G$. Then the group algebra $K[G]$ operates on the additive structure of $L$ in the obvious way and we define the associated order of the extension $L/K$ to be
$$\mathcal{A}_{L/K} := \{ x \in K[G] \, | \, x(\mathcal{O}_L) 
\subseteq \mathcal{O}_L \} ,$$
where $\mathcal{O}_{L}$ is the ring of integers of $L$.
A natural problem that arises in Galois module theory is that of determining the associated order $\mathcal{A}_{L/K}$ and the structure of $\mathcal{O}_{L}$ as an $\mathcal{A}_{L/K}$-module. \par

In the case where $K=\Q$ and $L/\Q$ is abelian, this problem was solved by Leopoldt 
(see \cite{leopoldt}; Lettl gives a simplified proof in \cite{lettl-global}). In particular, Leopoldt showed that $\mathcal{O}_{L}$ is a free 
$\mathcal{A}_{L/\Q}$-module (necessarily of rank 1) and gave explicit descriptions of $\mathcal{A}_{L/\Q}$ and of $\alpha \in \mathcal{O}_{L}$ such that $\mathcal{O}_{L}=\mathcal{A}_{L/\Q} \cdot \alpha$.

In the same situation ($K=\Q$ and $L/\Q$ abelian), one can ask whether $\mathcal{O}_{L}$ has a normal integral basis, i.e. whether it is free as a $\Z[G]$-module 
(again, necessarily of rank $1$). The Hilbert-Speiser Theorem says that this is the case if and only if $L$ is of square-free conductor or, equivalently, the extension $L/\Q$ is tamely ramified. In the more general setting, Noether showed that $\mathcal{A}_{L/K}=\mathcal{O}_{K}[G]$ precisely when the extension $L/K$ is at most tamely ramified (see \cite{noether}). Furthermore, $\mathcal{A}_{L/K}$ is the only order of $K[G]$ over which 
$\mathcal{O}_{L}$ can possibly be free of rank $1$. Hence Leopoldt's Theorem can be viewed as a generalization of the Hilbert-Speiser Theorem. \par

A number field $K$ is called a Hilbert-Speiser field if each at most tamely ramified finite abelian extension $L/K$ has a relative normal integral basis, i.e. $\mathcal{O}_{L}$ is a free $\mathcal{O}_{K}[G]$-module. In \cite{hilbert-speiser}, Greither, Replogle, Rubin and Srivastav show that $\Q$ is the only such number field. In this paper, we consider 
an analogous problem in which the restriction on ramification is relaxed and the requirement that $L/K$ is abelian is strengthened. 

\begin{definition}
An extension $L/K$ of absolutely abelian number fields (i.e. $L/\Q$ abelian) 
is \emph{Leopoldt} if $\mathcal{O}_{L}$ is a free 
$\mathcal{A}_{L/K}$-module. An abelian number field $K$ is \emph{Leopoldt} if every extension 
$L/K$ with $L/\Q$ abelian is Leopoldt in the sense above.
\end{definition}

In particular, Leopoldt's Theorem shows that $\Q$ is Leopoldt. However, unlike the case of Hilbert-Speiser fields, 
$\Q$ is not the only such field. A series of successively sharper results of Chan \& Lim, Bley, and Byott \& Lettl (see \cite{cl}, \cite{bley} and \cite{byott-lettl}, respectively) culminates in the proof that the $n$-th cyclotomic field $\Q^{(n)}$ is Leopoldt for every $n$. In \cite{gomez}, G\'omez Ayala shows that any quadratic extension $L/K$ with $L$ of 
$p$-power conductor over $\Q$ ($p$ an odd prime) is Leopoldt (also see \cite{lettl-quad}). In this paper, we build upon these results by proving the following theorems.

\begin{theorem}\label{maintheorem1}
Let $r,s \geq 1$ and $p$ be an odd prime.
Let $L/K$ be an extension with $[L:K]=p^{s}$ and $L \subseteq \Q^{(p^{r+s})}$.
Then $L/K$ is Leopoldt. 
\end{theorem}

\begin{theorem}\label{maintheorem2}
Let $r \geq 1$, $s \geq 0$ and $p$ be an odd prime. Then $\Q^{(p^{r+s})}/\Q^{(p^{r})+}$ is Leopoldt.
\end{theorem}

\begin{theorem}\label{maintheorem3}
Let $m \in \N$ with prime factorization $m=p_{1}^{r_{1}}p_{2}^{r_{2}} \cdots p_{k}^{r_{k}}$
where $p_{1}=2$ and $r_{1} \neq 1$. 
Let $K$ be the compositum of the collection of number fields 
$ \{ K_{i} \}_{i=1}^{k}$ where $K_{i} = \Q^{(p^{r_{i}})}$ or $\Q^{(p^{r_{i}})+}$ for $i \geq 2$ and
$K_{1}=\Q^{(2^{r_{1}})}$. Then $K$ is Leopoldt.
\end{theorem}

Note that these results can be made explicit (see Theorems \ref{cyclotomic-extension} and 
\ref{max-real-extension}). 
In other words, for a Leopoldt extension $L/K$ we can
give an explicit generator $\alpha \in \mathcal{O}_{L}$ such that 
$\mathcal{O}_{L}= \mathcal{A}_{L/K} \cdot \alpha$. 

Theorem \ref{maintheorem3} follows from Theorems \ref{maintheorem1} and \ref{maintheorem2}
(plus certain cases of the already known results mentioned above) in a relatively straightforward manner as shown
in Section \ref{piecing}. One important tool used in the proof of Theorems \ref{maintheorem1}
and \ref{maintheorem2}
is the result of \cite{lettl-local} in which Lettl shows that for any finite extension of $p$-adic fields 
$L_{\mathfrak{P}}/K_{\mathfrak{p}}$ where $L_{\mathfrak{P}}/\Q_{p}$ is abelian, 
$\mathcal{O}_{L_{\mathfrak{P}}}$ is a free $\mathcal{A}_{L_{\mathfrak{P}}/K_{\mathfrak{p}}}$-module. 
We then use Roiter's Lemma (see 31.6 of \cite{curtisandreiner}) to reduce the problem to the computation of
certain discriminants, which we do by explicitly calculating resolvents.

It is natural to ask whether every finite abelian extension of $\Q$ is Leopoldt. 
One can trivially deduce that the answer is no from results of Brinkhuis (\cite{brink}, \cite{brink-real}, \cite{brink-4}) 
or Cougnard (\cite{cougnard}), which give many examples of absolutely abelian tame extensions that 
have no relative normal integral basis. In analogy with the result that $\Q$ is the only Hilbert-Speiser field 
(\cite{hilbert-speiser}), we would ultimately like to give a classification of all Leopoldt fields. This is a problem 
that might be explored in a future paper.

\section{Piecing Extensions Together}\label{piecing}

We show that Theorem \ref{maintheorem3} follows from Theorems \ref{maintheorem1} 
and \ref{maintheorem2}, plus
Leopoldt's original theorem along with a special case ($2$-power conductor) 
of its previously known relative generalizations.

\begin{lemma}\label{bl-stick}
Let $L_{1}$ and $L_{2}$ be finite Galois extensions of a number field $K$. Put $L=L_{1}L_{2}$
and suppose that $L_{1}$ and $L_{2}$ are arithmetically disjoint over $K$ (i.e. 
$\mathcal{O}_{L} = \mathcal{O}_{L_{1}} \tensor_{\mathcal{O}_{K}} \mathcal{O}_{L_{2}}$,
or equivalently, $L_{1}$ and $L_{2}$ are linearly disjoint over $K$ and 
$(\disc(\mathcal{O}_{L_{1}}/\mathcal{O}_{K}),\disc(\mathcal{O}_{L_{2}}/\mathcal{O}_{K}))=\mathcal{O}_{K}$
- see III.2.13 of \cite{ft}). Then:
\begin{enumerate}
\item We have 
$\mathcal{A}_{L/L_{2}} = \mathcal{A}_{L_{1}/K} \tensor_{\mathcal{O}_{K}} \mathcal{O}_{L_{2}}$
and 
$\mathcal{A}_{L/K}= \mathcal{A}_{L_{1}/K} \tensor_{\mathcal{O}_{K}} \mathcal{A}_{L_{2}/K}$.
\item If there exists some $\alpha_{1} \in \mathcal{O}_{L_{1}}$ with 
$\mathcal{O}_{L_{1}}= \mathcal{A}_{L_{1}/K} \cdot \alpha_{1}$, then 
$\mathcal{O}_{L}=\mathcal{A}_{L/L_{2}} \cdot (\alpha_{1} \tensor 1)$. \\
If there also exists $\alpha_{2} \in \mathcal{O}_{L_{2}}$ with 
$\mathcal{O}_{L_{2}}= \mathcal{A}_{L_{2}/K} \cdot \alpha_{2}$, then
$\mathcal{O}_{L}=\mathcal{A}_{L/K} \cdot (\alpha_{1} \tensor \alpha_{2})$.
\end{enumerate}
\end{lemma}

\begin{proof}
This is Lemma 5 of \cite{byott-lettl}; the proof is immediate.
\end{proof}

\begin{lemma}\label{my-stick}
Let $L_{1}/K_{1}$ and $L_{2}/K_{2}$ both be Galois extensions of number fields.
Suppose that $L_{1}$ and $L_{2}$ are arithmetically disjoint over $\Q$ and that for $i=1,2$
there exist $\alpha_{i} \in \mathcal{O}_{L_{i}}$ such that 
$\mathcal{A}_{L_{i}/K_{i}} \cdot \alpha_{i} = \mathcal{O}_{L_{i}}$. Then
$\mathcal{O}_{L_{1}L_{2}}=\mathcal{A}_{L_{1}L_{2}/K_{1}K_{2}} \cdot (\alpha_{1} \tensor \alpha_{2})$
and  
$$
\mathcal{A}_{L_{1}L_{2}/K_{1}K_{2}} = 
\left( \mathcal{A}_{L_{1}/K_{1}} \tensor_{\mathcal{O}_{K_{1}}} \mathcal{O}_{K_{1}K_{2}} \right)
\tensor_{\mathcal{O}_{K_{1}K_{2}}} 
\left( \mathcal{A}_{L_{2}/K_{2}} \tensor_{\mathcal{O}_{K_{2}}} \mathcal{O}_{K_{1}K_{2}} \right)
.$$
\end{lemma}

\begin{proof}
We give a slightly modified version of the proof dealing with the cyclotomic case given in
Section 5 of \cite{cl}. Consider the following field diagram.
$$
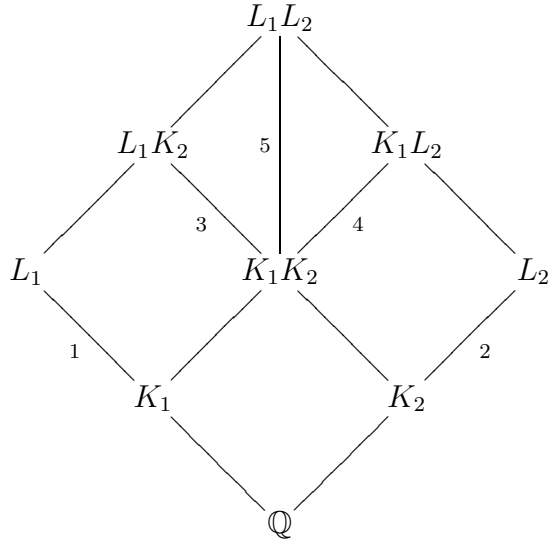
\xymatrix@1@!0@=48pt {
& & L_{1}L_{2} \ar@{-}[dr] \ar@{-}[dl] \ar@{-}[dd]_{5} & &  \\
& L_{1}K_{2}  \ar@{-}[dl] \ar@{-}[dr]_{3} & & K_{1}L_{2} \ar@{-}[dr] \ar@{-}[dl]^{4} & \\
L_{1}  \ar@{-}[dr]_{1} & & K_{1}K_{2} \ar@{-}[dl] \ar@{-}[dr] & & L_{2} \ar@{-}[dl]^{2} \\
& K_{1} \ar@{-}[dr] & & K_{2} \ar@{-}[dl]\\
& & \Q & &
}
$$
The hypothesis that $L_{1}$ and $L_{2}$ are arithmetically disjoint over $\Q$ implies that the
relevant pairs of fields are arithmetically disjoint over the relevant base field throughout. \par
We use Lemma \ref{bl-stick} to move from extensions (1) and (2) to (3) and (4)
respectively, and again to arrive at the result for extension (5). 
\end{proof}

\begin{prop}\label{wildtrace}
Let $L/K$ be an extension of number fields. Then $\tr_{L/K}(\mathcal{O}_L)$ is an ideal of $\mathcal{O}_K$. Suppose further that $L/K$ is Galois and let $\mathfrak{p}$ be a (non-zero) prime of $\mathcal{O}_K$. Then $\mathfrak{p} \mid \tr_{L/K}(\mathcal{O}_L)$ if and only if $\mathfrak{p}$ is wildly ramified in $L/K$.
\end{prop}

\begin{proof}
See \cite{mj}. Alternatively, this follows from Lemma 2 in Section 5 of \cite{lf} and the fact that the extension of residue fields in question is separable.
\end{proof}

\begin{corollary}\label{wildtracecor}
If $L/K$ is a Galois extension of number fields, then $L/K$ is at most tamely ramified
if and only if $\tr_{L/K}(\mathcal{O}_L)=\mathcal{O}_K$.
\end{corollary}

\begin{corollary}\label{tame-trace-down}
Let $M/L$ and $N/L$ be Galois extensions of number fields with
$L \subset M \subset N$ and $N/M$ at most tamely ramified. 
Put $G=\Gal(N/L)$ and $H=\Gal(M/L)$ and let $\pi:L[G] \rightarrow L[H]$ denote
the $L$-linear map induced by the natural projection $G \rightarrow H$. Suppose that
$\mathcal{O}_{N}=\mathcal{A}_{N/L}\cdot \alpha$ for some 
$\alpha \in \mathcal{O}_{N}$. 
Then $\mathcal{A}_{M/L} = \pi(\mathcal{A}_{N/L})$ 
and $\mathcal{O}_{M}=\mathcal{A}_{M/L}\cdot \tr_{N/M}(\alpha)$.
\end{corollary}

\begin{proof}
By Corollary \ref{wildtracecor},
 $\tr_{N/M}(\mathcal{O}_{N}) = \mathcal{O}_{M}$, and so
$$ \mathcal{O}_{M} = \tr_{N/M}(\mathcal{O}_{N})
= \tr_{N/M}(\mathcal{A}_{N/L} \cdot \alpha) 
= \mathcal{A}_{N/L}(\tr_{N/M}(\alpha))
= \pi(\mathcal{A}_{N/L}) \cdot \tr_{N/M}(\alpha).$$
(Note that the statement and proof are essentially the same as for Lemma 6 of \cite{byott-lettl}.)
\end{proof}

\begin{proof}[Proof of Theorem \ref{maintheorem3}]
Recall that $K$ is the compositum of the collection of number fields $\{ K_{i} \}_{i=1}^{k}$ where
$K_{i}=\Q^{(p^{r_{i}})}$ or $K_{i}=\Q^{(p^{r_{i}})+}$ for $i \geq 2$ and $K_{1}=\Q^{(2^{r_{1}})}$ 
($p_{1}=2$ and $r_{1} \neq 1$), and that we wish to show that $K$ is Leopoldt. In other words,  
given a finite extension $L/K$ with $L/\Q$ abelian, we need to find $\alpha \in \mathcal{O}_{L}$
such that $\mathcal{O}_{L}= \mathcal{A}_{L/K} \cdot \alpha$.

Let $n=p_{1}^{r_{1}+s_{1}}p_{2}^{r_{2}+s_{2}} \cdots p_{k}^{r_{k}+s_{k}}$ be 
the prime factorization of the conductor of $L$. Note that we may have $r_{i}=0$ for some $i$. For $1 \leq i \leq k$, let $L_{i}=\Q^{(p_{i}^{r_{i}+s_{i}})}$.
For each $1 \leq i \leq k$, there exists $\beta_{i} \in \mathcal{O}_{L_{i}}$
such that $\mathcal{O}_{L_{i}}=\mathcal{A}_{L_{i}/K_{i}} \cdot \beta_{i}$. We see this as follows: for 
$r_{i}=0$ (i.e. $K_{i}=\Q$), use Leopoldt's original theorem 
(see \cite{leopoldt} or the simplified proof in \cite{lettl-global});
for $i=1$ and $r_{1} \geq 2$ use a special case of the main result of \cite{cl} (or its generalizations in \cite{bley}
or \cite{byott-lettl}); and for all remaining cases, use Theorems \ref{maintheorem1} and \ref{maintheorem2}.
Successive applications of Lemma
\ref{my-stick} give $\mathcal{O}^{(n)}=\mathcal{A}_{\Q^{(n)}/K} \cdot \beta$ where 
$\beta= \beta_{1}\beta_{2} \cdots \beta_{k}$.
If $n$ is odd or if $\Q(i) \subseteq K$, then $\Q^{(n)}/L$ is at most tamely ramified
(see Proposition 3.7 and Remark 3.8 of \cite{me}) and so we have 
$\mathcal{O}_{L}=\mathcal{A}_{L/K} \cdot \tr_{\Q^{(n)}/L}(\beta)$
by Corollary \ref{tame-trace-down}.

Now suppose $n$ is even and that $\Q(i) \nsubseteq K$. Then $r_{1}=0$
(recall $p_{1}=2$). Assuming the notation and results of Sections 2 and 3 of \cite{me}, let $X$ be the group of Dirichlet characters associated to $L$, and let $L'$ be the field associated to $X_{2}$. Then by Leopoldt's original theorem, there exists $\gamma_{1} \in \mathcal{O}_{L'}$
such that $\mathcal{O}_{L'} = \mathcal{A}_{L'/\Q} \cdot \gamma_{1}$.
Let $t=p_{2}^{r_{2}+s_{2}} \cdots p_{k}^{r_{k}+s_{k}}$.
As in the odd case, there exists $\gamma_{2} \in \mathcal{O}^{(t)}$ such that
$\mathcal{O}^{(t)} = \mathcal{A}_{\Q^{(t)}/K} \cdot \gamma_{2}$. 
So by Lemma \ref{bl-stick} we have
$\mathcal{O}_{L'\Q^{(t)}}=\mathcal{A}_{L'\Q^{(t)}/K} \cdot \gamma_{1}\gamma_{2}$.
However, the extension $L'\Q^{(t)}/L$ is at most tamely ramified (for primes above odd rational primes, use Proposition 3.7 of \cite{me}; for primes above $2$, use Theorem 3.5 of \cite{wash}), and so by Corollary \ref{tame-trace-down}, we have 
$\mathcal{O}_{L}= \mathcal{A}_{L/K}\cdot \tr_{L'\Q^{(t)}/L}(\gamma_{1}\gamma_{2})$.
\end{proof}

\begin{remark}
Using the proof above and the results to which it refers, one can give explicit descriptions
of $\alpha$ and $\mathcal{A}_{L/K}$. Note, however, that there are in fact many possible 
choices for $\alpha$: if $\mathcal{O}_{L} = \mathcal{A}_{L/K} \cdot \alpha$, then
for any unit $u$ of $\mathcal{A}_{L/K}$, we have 
$\mathcal{O}_{L} = \mathcal{A}_{L/K} \cdot (u \cdot \alpha)$. 
\end{remark}

\begin{remark}
Theorem \ref{maintheorem3} can be easily extended using Lemma \ref{my-stick} to ``add on''
non-cyclotomic extensions of the form described in Theorem \ref{maintheorem1} (in the proof above, we only use the cyclotomic case of Theorem \ref{maintheorem1}) and quadratic extensions of the forms described 
in \cite{gomez} and \cite{lettl-quad}. However, the resulting statement would then be rather messy 
as it would be about certain Leopoldt extensions rather than Leopoldt fields.
\end{remark}

\section{From Local Freeness to Global Freeness}\label{local-global-freeness}

Throughout this section, let $L/K$ be a Galois extension of number fields with $G=\Gal(L/K)$ (not necessarily abelian). It is well-known that the ring of integers $\mathcal{O}_{L}$ being locally free over the associated order 
$\mathcal{A}_{L/K}$ is a necessary but not a sufficient condition for $\mathcal{O}_{L}$ to be free over 
$\mathcal{A}_{L/K}$. However, the proposition below says that in some sense local freeness is ``not too far'' 
from global freeness.

Let $[A:B]=[A:B]_{\mathcal{O}_{K}}$ denote 
the $\mathcal{O}_{K}$-module index of $B$ in $A$. (In fact, one could instead set 
$[A:B]=[A:B]_{\mathcal{O}_{F}}$ for any subfield $F$ of $K$.) For background material 
on module indices and discriminants, see Sections II.4 and III.2 of \cite{ft}.

\begin{prop}\label{special-roiter}
Suppose $\mathcal{O}_{L}$ is locally free over $\mathcal{A}_{L/K}$.
Then given any non-zero ideal $I$ of $\mathcal{O}_{K}$, there exists $\beta \in \mathcal{O}_{L}$ 
(depending on $I$) such that 
$I + [\mathcal{O}_{L}:\mathcal{A}_{L/K}\cdot \beta] = \mathcal{O}_{K}$.
\end{prop}

\begin{proof}
We use Roiter's Lemma as stated in 31.6 of \cite{curtisandreiner}
(alternatively, see Theorem 27.1 of \cite{reiner}). 
Let $R=\mathcal{O}_{K}$, $\Lambda=M=\mathcal{A}_{L/K}$ and $N=\mathcal{O}_{L}$.
Both $M$ and $N$ are $\Lambda$-lattices and are in the same genus
since $\mathcal{O}_{L}$ is locally free over $\mathcal{A}_{L/K}$.
Roiter's Lemma says that there exists a $\Lambda$-exact sequence
$$ 0 \rightarrow M \rightarrow N \rightarrow T \rightarrow 0$$
for some $\Lambda$-module $T$ such that $I+\mathrm{ann}_{R}(T)=R$. Let $S$ be
the image of $M$ in $N$ under the map of the above sequence. Then $S$ is a free
$\Lambda$-module of rank $1$ contained in $N=\mathcal{O}_{L}$ and so there exists
some $\beta \in S \subseteq N = \mathcal{O}_{L}$ such that 
$S=\Lambda \cdot \beta = \mathcal{A}_{L/K} \cdot \beta$. Hence we have an exact 
sequence
$$  0 \rightarrow \mathcal{A}_{L/K} \cdot \beta \rightarrow \mathcal{O}_{L} \rightarrow T \rightarrow 0 $$
where $I+\mathrm{ann}_{\mathcal{O}_{K}}(T)=\mathcal{O}_{K}$. 
The desired result now follows once one observes that
$$\mathrm{rad}([\mathcal{O}_{L}:\mathcal{A}_{L/K}\cdot \beta])
= \mathrm{rad}(\mathrm{ann}_{\mathcal{O}_{K}}(T)).$$
\end{proof}

The following proposition says that local freeness together with certain extra conditions implies global freeness.

\begin{prop}\label{sneaky-trick}
Suppose that $\mathcal{O}_{L}$ is locally free over $\mathcal{A}_{L/K}$.  
Let $\alpha \in \mathcal{O}_{L}$.
Then $\mathcal{O}_{L}=\mathcal{A}_{L/K} \cdot \alpha$ if and only if
$[\mathcal{O}_{L}:\mathcal{O}_{K}[G] \cdot \beta] \subseteq 
[\mathcal{O}_{L}:\mathcal{O}_{K}[G] \cdot \alpha]$ for all $\beta \in \mathcal{O}_{L}$.
\end{prop}

\begin{proof}
Let $S = \{ \beta \in \mathcal{O}_{L} \, | \, [\mathcal{O}_{L}:\mathcal{O}_{K}[G] \cdot \beta] \neq 0\}$
and observe that this is precisely the set of elements of $\mathcal{O}_{L}$ that generate a normal
basis for the field extension $L/K$. Note we can assume without loss of generality that $\alpha \in S$.
For every $\beta \in S$, the map
$$ \psi : K[G] \rightarrow L, \quad x \mapsto x \cdot \beta$$
is in fact an isomorphism (in particular, of $K$-vector spaces) which restricts to an 
isomorphism of $\mathcal{O}_{K}$-modules $\mathcal{A}_{L/K} \rightarrow \mathcal{A}_{L/K} \cdot \beta$. 
Thus for every $\beta \in S$, we have
$$ [\mathcal{O}_{L}:\mathcal{O}_{K}[G] \cdot \beta]
= [\mathcal{O}_{L}:\mathcal{A}_{L/K} \cdot \beta]
[\mathcal{A}_{L/K} \cdot \beta : \mathcal{O}_{K}[G] \cdot \beta] = 
[\mathcal{O}_{L}:\mathcal{A}_{L/K} \cdot \beta]
[\mathcal{A}_{L/K} : \mathcal{O}_{K}[G] ]. $$
In particular, this holds for $\beta=\alpha$.
Hence for all $\beta \in \mathcal{O}_{L}$ (including $\beta \notin S)$ 
$[\mathcal{O}_{L}:\mathcal{O}_{K}[G] \cdot \beta] \subseteq 
[\mathcal{O}_{L}:\mathcal{O}_{K}[G] \cdot \alpha]$
if and only if 
$[\mathcal{O}_{L}:\mathcal{A}_{L/K} \cdot \beta] \subseteq  [\mathcal{O}_{L} : \mathcal{A}_{L/K} \cdot \alpha]$.
However, $[\mathcal{O}_{L}:\mathcal{A}_{L/K} \cdot \beta] \subseteq  [\mathcal{O}_{L} : \mathcal{A}_{L/K} \cdot \alpha]$ for all $\beta \in \mathcal{O}_{L}$ if and only if $[\mathcal{O}_{L} : \mathcal{A}_{L/K} \cdot \alpha]=(1)$. 
To see this, note that one direction is trivial and that the other must hold in order to avoid a contradiction 
with Proposition \ref{special-roiter}. 
\end{proof}

We now rephrase Proposition \ref{sneaky-trick} in terms of discriminants.

\begin{prop}\label{disc-contain}
Suppose that $\mathcal{O}_{L}$ is locally free over $\mathcal{A}_{L/K}$.
Let $\alpha \in \mathcal{O}_{L}$. Then $\mathcal{O}_{L} = \mathcal{A}_{L/K} \cdot \alpha$
if and only if for every $\beta \in \mathcal{O}_{L}$, we have
$$ \disc( (\mathcal{O}_{K}[G] \cdot \beta) / \mathcal{O}_{K} ) \subseteq 
\disc( (\mathcal{O}_{K}[G] \cdot \alpha) / \mathcal{O}_{K} ).$$
\end{prop}

\begin{proof}
Using a well-known lemma relating discriminants to module indices (see III.2.4 of \cite{ft}), for all
$\beta \in \mathcal{O}_{L}$ (including $\beta=\alpha$) we have
$$ \disc( (\mathcal{O}_{K}[G] \cdot \beta) / \mathcal{O}_{K} ) 
= [\mathcal{O}_{L}:(\mathcal{O}_{K}[G] \cdot \beta)]^{2} \cdot \disc(\mathcal{O}_{L}/\mathcal{O}_{K}).$$
The desired result now follows from Proposition \ref{sneaky-trick}.
\end{proof}

\begin{prop}\label{shift-down}
Let $K$ be a number field with linearly disjoint extensions $L$ and $E$. Let $F=LE$
and suppose that $F/E$ (and so $L/K$) is Galois with Galois group $G$.
$$
\xymatrix@1@!0@=24pt { 
& & F \\
L \ar@{-}[urr] & & \\
& & E \ar@{-}[uu]_{G} \\
K \ar@{-}[uu]^{G} \ar@{-}[urr] & & \\
}
$$
Suppose that $\mathcal{O}_{L}$ and $\mathcal{O}_{F}$ are locally free over $\mathcal{A}_{L/K}$
and $\mathcal{A}_{F/E}$ respectively and that there exists $\alpha \in \mathcal{O}_{L}$ 
such that 
$\mathcal{O}_{F} = \mathcal{A}_{F/E} \cdot \alpha$. Then 
$\mathcal{O}_{L} = \mathcal{A}_{L/K} \cdot \alpha$.
\end{prop}

\begin{proof}
First note that the restriction map
$$ \mathrm{res} : G = \Gal(F/E) \rightarrow \Gal(L/K), \quad
\sigma \mapsto \sigma |_{L}$$
is an isomorphism by the linear disjointness hypothesis. 
For every $\beta \in \mathcal{O}_{L} \subseteq \mathcal{O}_{F}$, we have 
$$ \disc( (\mathcal{O}_{E}[G] \cdot \beta) / \mathcal{O}_{E} ) \subseteq 
\disc( (\mathcal{O}_{E}[G] \cdot \alpha) / \mathcal{O}_{E} ).$$
However, since 
$\mathcal{O}_{E}[G] = \mathcal{O}_{E} \tensor_{\mathcal{O}_{K}} \mathcal{O}_{K}[G]$
we have 
$$
\disc( (\mathcal{O}_{E}[G] \cdot \beta) / \mathcal{O}_{E})
= \mathcal{O}_{E} \tensor_{\mathcal{O}_{K}} \disc( (\mathcal{O}_{K}[G] \cdot \beta) / \mathcal{O}_{K})
$$
for all $\beta \in \mathcal{O}_{L}$ (including $\beta = \alpha$). It is now easy to see that 
$$ \disc( (\mathcal{O}_{K}[G] \cdot \beta) / \mathcal{O}_{K} ) \subseteq 
\disc( (\mathcal{O}_{K}[G] \cdot \alpha) / \mathcal{O}_{K} )$$
for every $\beta \in \mathcal{O}_{L}$ and so the desired result now follows from
Proposition \ref{disc-contain}.
\end{proof}

\section{Resolvents}\label{resolvents}

In this section, we define resolvents and use them to give a formula for certain discriminants.
We then derive a criterion in terms of resolvents for trivial Galois module structure of abelian extensions which 
will be a key element of the proofs of the main results.

\begin{definition}
For every $k \in \N$ choose a primitive $k$-th root of unity $\zeta_{k} \in \C$ in such a way that
for all $k,l \in \N$ with $k|l$ we have $\zeta_{l}^{l/k}=\zeta_{k}$. Let this choice be fixed for the
rest of the paper.
\end{definition}

Let $L/K$ be an abelian extension of number fields with Galois group $G$ and let $\widehat{G}$ be the
character group of $G$. Let $n$ be the exponent of $G$ and note that the elements of $\widehat{G}$ can be considered as group homomorphisms $\chi: G \rightarrow L(\zeta_{n})^{\times}$.

\begin{definition}\label{resolvent}
Let $\alpha \in L$ and $\chi \in \widehat{G}$. We define
$$ \langle \alpha \, | \, \chi \rangle = \langle \alpha \, | \, \chi \rangle_{L/K} 
:= \sum_{g \in G} \chi(g^{-1})g(\alpha) \in L(\zeta_{n})$$ 
to be the \emph{resolvent} attached to $\alpha$ and $\chi$.
\end{definition}

\begin{prop}\label{abelian-group-det}
Let $\alpha \in \mathcal{O}_{L}$. Then
$$ \disc((\mathcal{O}_{K}[G] \cdot \alpha) / \mathcal{O}_{K}) = \disc( \{ g(\alpha) : g \in G \} / \mathcal{O}_{K}) =\prod_{\chi \in \widehat{G}} \left( \langle \alpha \, | \, \chi \rangle \right)^{2}.$$
\end{prop}

\begin{proof}
Lemma 5.26 (a) in \cite{wash} gives
$ \det((hg^{-1}(\alpha))_{h,g \in G}) = \prod_{\chi \in \widehat{G}} \langle \alpha \, | \, \chi \rangle$
for $\alpha \in L$ (this is known as the abelian group determinant formula). Proposition 2.25 of \cite{milne} says that $\disc( \{ g(\alpha) : g \in G \} / \mathcal{O}_{K}) =  (\det((hg^{-1}(\alpha))_{h,g \in G}))^{2}$. Thus we obtain the desired result.
\end{proof}

\begin{corollary}\label{resolvent-divide}
Suppose $\mathcal{O}_{L}$ is locally free over $\mathcal{A}_{L/K}$. Let $\alpha \in \mathcal{O}_{L}$.
Then 
$$ \prod_{\chi \in \widehat{G}} \langle \alpha \, | \, \chi \rangle \textrm{ divides } 
\prod_{\chi \in \widehat{G}}  \langle \beta \, | \, \chi \rangle \quad \forall \beta \in \mathcal{O}_{L}$$
if and only if $\mathcal{O}_{L} = \mathcal{A}_{L/K} \cdot \alpha$. (Note: everything divides zero.)
\end{corollary}

\begin{proof}
Follows from Propositions \ref{disc-contain} and \ref{abelian-group-det}.
\end{proof}

We now give a reduction step useful for computing resolvents.

\begin{lemma}\label{res-trace}
Let $M$ be an intermediate field of $L/K$, $H = \Gal(L/M)$ and $\Delta=\Gal(M/K)$.
Suppose $\chi$ is a character of $G$ that is trivial on $H$. Then by abuse of notation we can consider 
$\chi$ as a character on $\Delta$, and for $\alpha \in L$ we have
$$  \langle \alpha \, | \, \chi \rangle_{L/K} =  \langle \tr_{L/M}(\alpha) \, | \, \chi \rangle_{M/K}.$$
\end{lemma}

\begin{proof}
Straightforward.
\end{proof}

\section{Sums of Roots of Unity}\label{roots-of-unity}

We compute certain sums of roots of unity, allowing us to calculate resolvents later on.

\begin{definition}
For any prime $p$, let $v_{p}(x)$ denote the $p$-adic valuation of $x \in \Z$.
\end{definition}

\begin{lemma}\label{root-sum}
Let $p$ be a prime, $a,b,k \in \Z$ and $k \geq 1$. 
Then
$$  \sum_{t=0}^{p^{k}-1} \zeta_{p^{k}}^{a+bt} = 
\left\{ 
\begin{array}{ll} 
\zeta_{p^{k}}^{a}p^{k} & \textrm{ if } v_{p}(b) \geq k, \\
0 & \textrm{ otherwise. }
\end{array} \right.
$$
\end{lemma}

\begin{proof}
Straightforward.
\end{proof}

\begin{lemma}\label{trace-extract}
Let $p$ be a prime, $n,m \geq 1$ ($n \geq 2$ if $p=2$) and $0 \leq k \leq n+m$.
Then $$ \tr_{\Q^{(p^{n+m})}/\Q^{(p^{n})}}(\zeta_{p^{k}}) = \left\{ 
\begin{array}{ll} 
\zeta_{p^{k}} p^{m} & \textrm{if } k \leq n, \\
0 & \textrm{otherwise.}
\end{array} \right.
$$
\end{lemma}

\begin{proof}
We have
$$ \tr_{\Q^{(p^{n+m})}/\Q^{(p^{n})}}(\zeta_{p^{k}}) = 
\sum_{i=0}^{p^{m}-1} \sigma_{1+ip^{n}}(\zeta_{p^{k}})
= \zeta_{p^{k}} \! \! \sum_{i=0}^{p^{m}-1} \zeta_{p^{k}}^{ip^{n}} $$
where $\sigma_{i}$ is the automorphism of $\Q(\zeta_{p^{n+m}})=\Q^{(p^{n+m})}$ 
mapping $\zeta_{p^{n+m}} \mapsto \zeta_{p^{n+m}}^{i}$ (and so $\zeta_{p^{k}} \mapsto \zeta_{p^{k}}^{i}$). 
The result now follows from Lemma \ref{root-sum}.
\end{proof}

\begin{theorem}[Gau\ss]\label{square-sum}
Let $n \in \N$. Let $\zeta=e^{2\pi i /n}$. Then
$$  \sum_{t=0}^{n-1} \zeta^{t^{2}} = 
\left\{ 
\begin{array}{rl} 
\sqrt{n} & \textrm{ if } n \equiv 1 \, \mod \, 4, \\
i\sqrt{n} & \textrm{ if } n \equiv 3 \, \mod \, 4.
\end{array} \right.
$$
\end{theorem}

\begin{proof}
See Theorem 99 of \cite{nag}.
\end{proof}

\begin{corollary}\label{top-level-zero}
Let $p$ be an odd prime, $k \geq 2$ and $\zeta=\zeta_{p^{k}}$. Then
$ \sum_{\substack{0 \leq t \leq p^{k}-1 \\ (p,t)=1}} \zeta^{t^{2}} = 0$.
\end{corollary}

\begin{proof}
Observe that 
$$ \sum_{t=0}^{p^{k}-1} \zeta^{t^{2}}  
= \sum_{\substack{0 \leq t \leq p^{k}-1 \\ (p,t)=1}} \zeta^{t^{2}}
+ \sum_{t=0}^{p^{k-1}-1} \zeta^{(pt)^{2}} 
= \sum_{\substack{0 \leq t \leq p^{k}-1 \\ (p,t)=1}} \zeta^{t^{2}}
+ p\sum_{t=0}^{p^{k-2}-1} \zeta^{p^{2}t^{2}}.
$$
The result for now follows by using Theorem \ref{square-sum} with $n=p^{k}$ and $n=p^{k-2}$.
(Strictly speaking, we have shown the result for $\zeta=e^{2\pi i /p^{k}}$, though it is clear that
the result still holds for any choice of $\zeta=\zeta_{p^{k}}$.)
\end{proof}

\begin{prop}\label{hensel}
Let $F(t) \in \Z_{p}[t]$, and let $\alpha_{0}$ be a simple root of $F(t) \, \mod \, p$.
Then there exists a root $\alpha$ of $F(t)$ with $\alpha \equiv \alpha_{0} \, \mod \, p$.  
\end{prop}

\begin{proof}
This is a corollary of Hensel's Lemma - see II.3.25.a of \cite{ft}.
\end{proof}

\begin{prop}
Let $p$ be an odd prime and let $a,b,k \in \Z$ with $k \geq 1$ and $(b,p)=1$.
Let $\zeta=\zeta_{p^{k}}$ and let $f \in (\Z/p^{k}\Z)[t]$. Define
$g:\Z/p^{k}\Z \rightarrow \Z/p^{k}\Z$ by $t \mapsto a + bt + pf(t)$.
Then $\sum_{t=0}^{p^{k}-1} \zeta^{g(t)} = \sum_{t=0}^{p^{k}-1} \zeta^{a+bt}= 0.$
\end{prop}

\begin{proof}
Since $(b,p)=1$, we see that $g$ is surjective $\mod \, p$. Hence
Proposition \ref{hensel} shows that $g$ is surjective and therefore bijective.
However,
$h:\Z/p^{k}\Z \rightarrow \Z/p^{k}\Z, \, t \mapsto a+bt$
is also bijective, and so the desired result now follows from Lemma \ref{root-sum}.
\end{proof}

\begin{corollary}\label{linear-zero}
Let $p$ be an odd prime and let $a,b,c,d,e,k \in \Z$ with $d,e,k \geq 1$ and $(b,p)=1$.
Let $\zeta=\zeta_{p^{k}}$ and define 
$g:\Z/p^{k}\Z \rightarrow \Z/p^{k}\Z$ by $t \mapsto 
a + bt + cp \sum_{j=2}^{d}\binom{t}{j}p^{(j-2)e} $.
Then 
$\sum_{t=0}^{p^{f}-1} \zeta^{g(t)} = 0$.
\end{corollary}

\begin{proof}
The main point is to check that $f(t)=\sum_{j=2}^{d}\binom{t}{j}p^{(j-2)e}$ is indeed an element
of $(\Z/p^{k}\Z)[t]$. One has to consider $p$-denominators of the binomial coefficients and 
a straightforward induction argument gives the desired result.
\end{proof}

\begin{prop}
Let $p$ be an odd prime and let $a,b,c,k \in \Z$ with $k \geq 0$ and $(c,p)=1$.
Let $\zeta=\zeta_{p^{k}}$ and let $f \in (\Z/p^{k}\Z)[t]$. Define
$g:\Z/p^{k}\Z \rightarrow \Z/p^{k}\Z$ by $t \mapsto a + bt + ct^{2} + pf(t)$.
Then 
$\Norm_{\Q^{(p^{k})}/\Q} \left( \sum_{t=0}^{p^{k}-1} \zeta^{g(t)} \right)
= \Norm_{\Q^{(p^{k})}/\Q} \left( \sum_{t=0}^{p^{k}-1} \zeta^{t^{2}} \right)
= p^{\frac{1}{2}k[\Q^{(p^{k})}:\Q]}
= p^{\frac{1}{2}k(p-1)p^{k-1}}$.
\end{prop}

\begin{proof}
Define $h$ to be the first three terms of $g$, that is $h(t)=a+bt+ct^{2}$.
Since multiplying $g(t)$ by the multiplicative inverse of $c$ in $\Z/p^{k}\Z$ 
(recall $(c,p)=1$) does not affect the norm of $\sum \zeta^{g(t)}$, we
can without loss of generality suppose that $c \equiv 1 \, \mod \, p^{k}$. Thus we have
$h(t) = a + bt + t^{2}$. Now completing the square, we see that there exists a change
of variables $x=t-b/2$ such that $h(x)=x^{2}$. (Note that without loss of generality
we can take the constant term to be zero, as again this does not affect the norm.)
Thus we have $g(x) = x^{2} + pf(x+b/2)$, and since $f$ was arbitrary, we can suppose that
in fact $g(x) = x^{2}+pf(x)$. 

If $k=0$ or $1$, the result now follows from Theorem \ref{square-sum} and the fact that $g = h$ in this case
(both are defined as functions on $\Z/p^{k}\Z$). We suppose that $k \geq 2$ and proceed by induction.

Let $H_1 := \{ j \in (\Z/p^{k}\Z)^{\times} \, | \, 1 \leq (j \, \mod \, p) \leq (p-1)/2 \}$, let 
$H_2 := (\Z/p^k\Z)^{\times}-H_1$ and let 
$H_3 := (p\Z/p^k\Z)$. Note that $H_1$, $H_2$ and $H_3$ form a partition
of $(\Z/p^k\Z)$, and that the congruence class $\mod \, p$ of an element 
$x \in (\Z/p^k\Z)$ uniquely determines the $H_i$ in which it lies.
Furthermore, $H_{1}=-H_{2}$ and so both the restrictions $h |_{H_1}$ and $h |_{H_2}$ 
are bijections onto their common image, $(\Z/p^k\Z)^{\times 2}$ (i.e. the squares in $(\Z/p^k\Z)^{\times}$).

Since $g \equiv h \, \mod \, p$ and $g'(x) \equiv 2x \not \equiv 0 \, \mod \, p$ for all 
$x \in (\Z/p^{k}\Z)^{\times}$, Proposition \ref{hensel} shows that $g(H_{1})=g(H_{2})=(\Z/p^f\Z)^{\times 2}=h(H_{1})=h(H_{2})$. Therefore for $i=1,2$ we have permutations $s_i:H_i \rightarrow H_i$ such that 
$h|_{H_i}(x) = g|_{H_i}(s_i(x))$. We can ``glue'' these together to give a permutation 
$s : (\Z/p^{k}\Z)^{\times} \rightarrow (\Z/p^{k}\Z)^{\times}$ such that $h(x)=g(s(x))$ for all 
$x \in (\Z/p^{k}\Z)^{\times}$. By Corollary \ref{top-level-zero}, we now have
$$ \sum_{\substack{0 \leq x \leq p^{k}-1 \\ (p,x)=1}} \zeta^{g(x)} = 
\sum_{\substack{0 \leq x \leq p^{k}-1 \\ (p,x)=1}} \zeta^{x^{2}} = 0
\quad \textrm{ and so } \quad \sum_{x=0}^{p^{k}-1} \zeta^{g(x)}  
= \sum_{x=0}^{p^{k-1}-1} \zeta^{g(px)}.$$
Write $f(x) = d + ex + x^{2}\tilde{f}(x)$ and let $\tilde{g}(x)=x^{2} +ex + px^{2}\tilde{f}(px)$.
Then we have
$$g(px)= p^{2}x^{2} + pf(px) = p^{2}x^{2} +pd + p^{2}ex + p^{3}x^{2}\tilde{f}(px)
 = pd+p^{2}\tilde{g}(x).$$
As before, we can suppose without loss of generality that $d=0$, and so
$$ \sum_{x=0}^{p^{k}-1} \zeta^{g(x)}  
= \sum_{x=0}^{p^{k-1}-1} \zeta^{g(px)}
=  \sum_{x=0}^{p^{k-1}-1} \zeta^{p^{2}\tilde{g}(x)}
= p \sum_{x=0}^{p^{k-2}-1} \zeta_{p^{k-2}}^{\tilde{g}(x)}.$$
Since $\tilde{g}$ is of the correct form and
$$ \sum_{x=0}^{p^{k}-1} \zeta^{x^{2}}  = \sum_{x=0}^{p^{k-1}-1} \zeta^{(px)^{2}}
=  p \sum_{x=0}^{p^{k-2}-1} \zeta_{p^{k-2}}^{x^{2}},$$
the result now follows by induction.
\end{proof}

\begin{corollary}\label{quadratic-non-zero}
Let $p$ be an odd prime and let $a,b,c,d,e,k \in \Z$ with $d \geq 2$, $e,k \geq 1$, 
$(c,p)=1$ (if $p=3$, assume $e \geq 2$).
Let $\zeta=\zeta_{p^{k}}$ and define 
$$
\textstyle{
g:\Z/p^{k}\Z \rightarrow \Z/p^{k}\Z \quad \textrm{by} \quad 
t \mapsto a + bt + c \sum_{j=2}^{d}\binom{t}{j}p^{(j-2)e}}.$$
Then 
$\Norm_{\Q^{(p^{k})}/\Q} \left( \sum_{t=0}^{p^{k}-1} \zeta^{g(t)} \right)
= p^{\frac{1}{2}k(p-1)p^{k-1}}$.
\end{corollary}

\begin{proof}
Write $$\textstyle{g(t) = a + bt + c \binom{t}{2} + c \sum_{j=3}^{d}\binom{t}{j}p^{(j-2)e} 
= a + (b-\frac{c}{2})t + \frac{c}{2}t^{2} +  cp \sum_{j=3}^{d}\binom{t}{j}p^{(j-2)e-1}.}$$  
The main point is to check that $f(t)=\sum_{j=3}^{d}\binom{t}{j}p^{(j-2)e-1}$ is indeed an element
of $(\Z/p^{k}\Z)[t]$. One has to consider $p$-denominators of the binomial coefficients and 
a straightforward induction argument gives the desired result except when 
$e=1$ and $p=3$.
\end{proof}

\section{Absolutely Abelian Extensions}

We recall two results for absolutely abelian extensions of number fields. 
The first of these results says that all absolutely abelian extensions of number fields
have the local freeness property which allows the results of Sections \ref{local-global-freeness} 
and \ref{resolvents} to be applied.

\begin{theorem}[\cite{lettl-local}]
Let $L_{\mathfrak{P}}/K_{\mathfrak{p}}$ be a finite extension of $p$-adic fields 
with $L_{\mathfrak{P}}/\Q_{p}$ abelian. Then 
$\mathcal{O}_{L_{\mathfrak{P}}}$ is a free 
$\mathcal{A}_{L_{\mathfrak{P}}/K_{\mathfrak{p}}}$-module of rank $1$.
\end{theorem}

\begin{corollary}[\cite{lettl-local}]\label{locally-free}
If $L/K$ is an extension of number fields with $L/\Q$ abelian, then
$\mathcal{O}_{L}$ is locally free over $\mathcal{A}_{L/K}$, i.e. for any prime
$\mathfrak{p}$ of $\mathcal{O}_{K}$,
$\mathcal{O}_{L, \mathfrak{p}}$ is a free $\mathcal{A}_{L/K, \mathfrak{p}}$-module
of rank $1$.
\end{corollary}

Let $L/K$ be an extension of number fields with $L/\Q$ abelian and let $G=\Gal(L/K)$.
The second result asserts that the associated order $\mathcal{A}_{L/K}$ is equal to $\mathcal{M}_{K[G]}$, 
the maximal $\mathcal{O}_{K}$-order of $K[G]$, in many cases of interest. Note that
$\mathcal{M}_{K[G]}$ is unique in this case because $G$ is abelian.

\begin{definition}
We say that an extension $N/M$ is \emph{totally wildly ramified} if each intermediate field
different from $M$ is wildly ramified above $M$.
\end{definition}

\begin{theorem}[\cite{byott-lettl}]\label{bl-max}
Let $L/K$ be a cyclic and totally wildly ramified extension of number fields with $L/\Q$ abelian. Suppose $L/K$ is linearly disjoint from the extension $\Q^{(m)}/K$, where $m$ 
denotes the conductor of $K$.
Then $\mathcal{A}_{L/K}=\mathcal{M}_{K[G]}$ where $G=\Gal(L/K)$.
\end{theorem}

\section{The Induction Step}\label{step}

We now specialize to the following situation:
\begin{itemize}
\item Let $r,s \geq 1$ and $p$ be an odd prime.
\item Let $V = \{ x \in \Z : v_{p}(x) = 0\}$.
\item Let $\zeta=\zeta_{p^{r+s}}$.
\item Let $L=\Q^{(p^{r+s})}$, $M=\Q^{(p^{r+s-1})}$ and $K=\Q^{(p^{r})}$ or $\Q^{(p^{r})+}$.
\item Let $G = \Gal(L/K)$, $H=\Gal(L/M)$ and $\Delta = \Gal(M/K)$.
\item Let $A$ be the subgroup of characters of $G$ that are trivial on $H$. Note that $A \cong \widehat{\Delta}$.
\item Let $B = \widehat{G} -A$. Note that this is the set of characters of $G$ of order $p^{s}$ or $2p^{s}$. 
\end{itemize}

\begin{remark}\label{exceptional}
We shall assume for the rest of the paper 
that $r \geq \min\{2,s\}$ if $p=3$ due to the exceptional case of Corollary \ref{quadratic-non-zero}
(note that this proposition is only used when $s>r$). However, we already know that extensions of the form 
$\Q^{(3^{1+s})}/\Q^{(3)}$ are Leopoldt by the main result of \cite{cl} (or its generalizations in \cite{bley} 
or \cite{byott-lettl}). Furthermore, the only
proper subfield of $\Q^{(3)}$ is $\Q^{(3)+}=\Q$, so results analogous to Corollary \ref{p-power-shifted} 
and Theorem \ref{max-real-extension} in the case $p=3$ and $r=1$ follow from Leopoldt's
original theorem \cite{leopoldt} (see \cite{lettl-global} for a simplified proof). Thus Theorems
\ref{maintheorem1} and \ref{maintheorem2} as stated still hold in this case.
\end{remark}

\begin{lemma}
Let $\sigma_{i}$ denote the field automorphism of $L=\Q(\zeta)$ defined
by $\zeta \mapsto \zeta^{i}$. Then
$$ G=\Gal(L/K)=\left\{ 
\begin{array}{ll} 
\{ \sigma_{1+ip^{r}} \, : \, 0 \leq i \leq p^{s}-1  \} = \langle \sigma_{1+p^{r}} \rangle 
& \textrm{ if } K=\Q^{(p^{r})}, \\
\{ \sigma_{1+ip^{r}} \, : \, 0 \leq i \leq p^{s}-1  \} \times \{ \sigma_1, \sigma_{-1} \} 
= \langle \sigma_{-1} \sigma_{1+p^{r}} \rangle
& \textrm{ if } K=\Q^{(p^{r})+}.
\end{array} \right.$$
\end{lemma}

\begin{proof}
Straightforward.
\end{proof}

\begin{lemma}\label{2p-split}
Let $b \in \Z$ and $\chi \in \widehat{G}$. Then
$$ \langle \zeta^{b} \, | \, \chi \rangle_{L/\Q^{(p^{r})+}} = \left\{ \begin{array}{ll}
\langle \zeta^{b} \, | \, \chi \rangle_{L/\Q^{(p^{r})}} - \langle \zeta^{-b} \, | \,  \chi  \rangle_{L/\Q^{(p^{r})}} & \textrm{ if }
\chi \textrm{ is of even order,} \\
\langle \zeta^{b} \, | \, \chi \rangle_{L/\Q^{(p^{r})}}  + \langle \zeta^{-b} \, | \, \chi \rangle_{L/\Q^{(p^{r})}}  & \textrm{ if }
\chi \textrm{ is of odd order.} \\
\end{array} \right. $$
\end{lemma}

\begin{proof}
If $\chi$ is of even order, we have
\begin{eqnarray*}
\langle \zeta^{b}  \, | \, \chi \rangle_{L/\Q^{(p^{r})+}} &=& \sum_{k=0}^{1} \sum_{j=0}^{p^{s}-1} 
\chi(((\sigma_{-1})^{k}(\sigma_{1+p^{r}})^{j})^{-1}) \sigma_{(-1)^{k}(1+p^{r})^{j}}(\zeta^{b}) \\
&=& \sum_{j=0}^{p^{s}-1} 
\chi((\sigma_{1+p^{r}})^{-j}) \sigma_{(1+p^{r})^{j}}(\zeta^{b})
- \sum_{j=0}^{p^{s}-1} 
\chi((\sigma_{1+p^{r}})^{-j}) \sigma_{(1+p^{r})^{j}}(\zeta^{-b}) \\
&=& \langle \zeta^{b} \, | \,  \chi \rangle_{L/\Q^{(p^{r})}} - \langle \zeta^{-b} \, | \,  \chi \rangle_{L/\Q^{(p^{r})}}.
\end{eqnarray*}
The calculation is almost identical for $\chi$ of odd order.
\end{proof}

\begin{lemma}\label{reduce-to-level-zero}
Let $\chi \in B$ and let $b \in \Z$ with $v_{p}(b)>0$. Then $\langle \zeta^{b} \, | \, \chi \rangle_{L/K}=0$.
\end{lemma}

\begin{proof}
Since $\chi \in B$, we have $\chi(\sigma_{1+p^{r}})=\zeta^{ap^{r}}$ for some $a \in V$.
Then using Definition \ref{resolvent} and noting the explanations below, we have
\begin{eqnarray*}
\langle \zeta^{b} \, | \, \chi \rangle_{L/\Q^{(p^{r})}}
& = & \sum_{j=0}^{p^{s}-1} \chi((\sigma_{1+p^{r}})^{-j})
(\sigma_{1+p^{r}})^{j}(\zeta^{b})
= \sum_{j=0}^{p^{s}-1} \zeta^{-ajp^{r}} 
\zeta^{b \sum_{i=0}^{j} \binom{j}{i}p^{ir}} \\
&=&  \sum_{j=0}^{p^{s}-1} \zeta^{-ajp^{r}} 
\zeta^{b \sum_{i=0}^{\lfloor s/r \rfloor +1} \binom{j}{i}p^{ir}} 
= \zeta^{b} \sum_{j=0}^{p^{s}-1} \zeta^{p^{r}\left((b-a)j+bp^{r}
\sum_{i=2}^{\lfloor s/r \rfloor +1} \binom{j}{i} p^{(i-2)r}\right)} = 0,
\end{eqnarray*}
where the last equality follows from Corollary \ref{linear-zero}, noting that $v_{p}(b-a)=0$. 
For the third equality above, note that $\binom{j}{i}=0$ for $i > j$ and that $p^{ir} \equiv 0 \, \mod \, p^{r+s}$ for 
$i \geq \lfloor s/r \rfloor + 2$. The $K=\Q^{(p^{r})+}$ case now follows from Lemma \ref{2p-split}.
\end{proof}

\begin{prop}\label{check-generator}
Suppose there exist $\gamma \in \mathcal{O}_{M}$, $f \in \N$ and $S \subseteq V$ such that 
\begin{enumerate}
\item $\mathcal{O}_{M} = \mathcal{A}_{M/K} \cdot \gamma$;
\item for every $\chi \in B$ and every $b \in V$, we have
$\Norm_{L/\Q}((\langle \zeta^{b} \, | \, \chi \rangle_{L/K})) = 0$ or $p^{f}$; and
\item for every $\chi \in B$, there exists a unique element
$b = b_{\chi} \in S$ such that $\langle \zeta^{b} \, | \, \chi \rangle_{L/K} \neq 0$.
\end{enumerate}
Then $\mathcal{O}_{L} = \mathcal{A}_{L/K} \cdot \alpha$ where $\alpha := \gamma + \sum_{b \in S} \zeta^{b}$.
\end{prop}

\begin{proof}
Let $\beta \in \mathcal{O}_{L}$. By Lemmas \ref{res-trace} and \ref{trace-extract}, we have
\begin{equation}\label{alpha-A-break}
\prod_{\chi \in A} \langle \alpha \, | \, \chi \rangle_{L/K}
= \prod_{\chi \in \widehat{\Delta}} \langle \tr_{L/M}(\alpha) \, | \, \chi \rangle_{M/K} 
= \prod_{\chi \in \widehat{\Delta}} \langle p \gamma \, | \, \chi \rangle_{M/K}
= p^{|\widehat{\Delta}|} \prod_{\chi \in \widehat{\Delta}} \langle \gamma \, | \, \chi \rangle_{M/K}.
\end{equation}
Let $\theta = p^{-1}\tr_{L/M}(\beta)$. Then $\theta \in \mathcal{O}_{M}$ by Lemma \ref{trace-extract}.
Furthermore, as above, we have
\begin{equation}\label{beta-A-break}
\prod_{\chi \in A} \langle \beta \, | \, \chi \rangle_{L/K} = 
p^{|\widehat{\Delta}|} \prod_{\chi \in \widehat{\Delta}} \langle \theta \, | \, \chi \rangle_{M/K}.
\end{equation}
However, since $\mathcal{O}_{M}=\mathcal{A}_{M/K} \cdot \gamma$, 
Corollary \ref{resolvent-divide} shows that (\ref{alpha-A-break}) divides 
(\ref{beta-A-break}).

Since $\gamma \in \mathcal{O}_{M}$, we can write $\gamma = \sum_{t \in T} \zeta^{pt}$ for some $T \subset \Z$.
Let $\mathfrak{p}=(1-\zeta)$ be the unique prime ideal in $\mathcal{O}_{L}=\Z[\zeta]$ above $p$.
Note that $\Norm_{L/\Q}(x)=p^{f}$ is equivalent to $(x)=\mathfrak{p}^{f}$ for any $x \in \mathcal{O}_{L}$.
Therefore by Lemma \ref{reduce-to-level-zero} and hypotheses (b) and (c), we have
\begin{eqnarray}\label{alpha-B-break}
\prod_{\chi \in B} \langle \alpha \, | \, \chi \rangle_{L/K}
&=& \prod_{\chi \in B} \left\langle 
{\textstyle \sum_{t \in T} \zeta^{pt} + \sum_{b \in S} \zeta^{b} } \, | \, \chi \right\rangle_{L/K}  \\
&=& \prod_{\chi \in B} \left\langle {\textstyle \sum_{b \in S} \zeta^{b} } \, | \, \chi \right\rangle_{L/K} 
\textrm{ generates } \frak{p}^{f|B|}. \nonumber
\end{eqnarray}
Similarly, we can write $\beta = \sum_{t \in T'} \zeta^{pt} + \sum_{s \in S'} \zeta^{s}$
for some $S' \subseteq V$ and $T' \subseteq \Z$. Then 
\begin{equation}\label{beta-B-break}
\prod_{\chi \in B} \langle \beta \, | \, \chi \rangle_{L/K}
= \prod_{\chi \in B} \left\langle 
{\textstyle \sum_{t \in T'} \zeta^{pt} + \sum_{s \in S'} \zeta^{s} } \, | \, \chi \right\rangle_{L/K} 
= \prod_{\chi \in B} \left\langle {\textstyle \sum_{s \in S'} \zeta^{s} } \, | \, \chi \right\rangle_{L/K} 
\in \frak{p}^{f|B|}.
\end{equation}
Therefore (\ref{alpha-B-break}) divides (\ref{beta-B-break}) and so 
$\prod_{\chi \in \widehat{G}} \langle \alpha \, | \, \chi \rangle_{L/K}$ divides 
$\prod_{\chi \in \widehat{G}} \langle \beta \, | \, \chi \rangle_{L/K}$
(note that $\widehat{G}$ is the disjoint union of $A$ and $B$). Since $\beta \in \mathcal{O}_{L}$
was arbitrary, the desired result now follows from Corollary \ref{resolvent-divide}.
\end{proof}

\section{Proof of Theorem \ref{maintheorem1}}\label{p-power}

We now prove an explicit version of Theorem \ref{maintheorem1}. We adopt the
notation given at the beginning of Section \ref{step}, taking $K=\Q^{(p^{r})}$.
We abbreviate $\langle - \, | \,  \chi \rangle_{L/K}$ to $\langle - \, | \,  \chi \rangle$. 

\begin{lemma}\label{pick}
Let $\chi \in B$ and define $a \in \N$ by $\chi (\sigma_{1+p^{r}}) = \zeta_{p^{s}}^{a}= \zeta^{ap^{r}}$.
If $p=3$, assume that $r \geq \min\{2,s\}$.
Let $b \in V$. Then 
$$
\Norm_{\Q^{(p^{r+s})}/\Q}( ( \langle  \zeta^{b} \, | \, \chi \rangle) ) = \left\{
\begin{array}{ll}
p^{s(p-1)p^{r+s-1}} & \textrm{if } s \leq r \textrm{ and } v_p(b-a) \geq s, \\
p^{\frac{1}{2}(s+r)(p-1)p^{r+s-1}} & \textrm{if } s > r \textrm{ and } v_p(b-a) \geq r, \\
0 & \textrm{otherwise.}
\end{array} \right.
$$
\end{lemma}

\begin{proof}
Exactly the same computation as in the proof of Lemma \ref{reduce-to-level-zero} gives 
$$\langle \zeta^{b} \, | \, \chi \rangle
= \zeta^{b} \sum_{j=0}^{p^{s}-1} \zeta^{p^{r}\left((b-a)j+bp^{r}
\sum_{i=2}^{\lfloor s/r \rfloor +1} \binom{j}{i} p^{(i-2)r}\right)}.$$
Suppose $s \leq r$. Then using Lemma \ref{root-sum} for the second equality, we have
$$
\langle \zeta^{b} \, | \, \chi \rangle = \zeta^{b} \sum_{j=0}^{p^{s}-1} \zeta^{p^{r}(b-a)j}
= \left\{
\begin{array}{ll}
\zeta^{b}p^{s} & \textrm{if } v_p(b-a) \geq s,\\
0 & \textrm{otherwise.}
\end{array} \right.
$$
Suppose $s > r$.
Let $l=v_{p}(b-a)$ and write $(b-a)=up^{l}$ where $(u,p)=1$. If $l<r$, we have 
$$\langle \zeta^{b} \, | \, \chi \rangle
= \zeta^{b} \sum_{j=0}^{p^{s}-1} \zeta^{p^{r+l}\left((uj+bp^{r-l}
\sum_{i=2}^{\lfloor s/r \rfloor +1} \binom{j}{i} p^{(i-2)r}\right)}
= \zeta^{b} p^{l} \sum_{j=0}^{p^{s-l}-1} \zeta^{p^{r+l}\left((uj+bp^{r-l}
\sum_{i=2}^{\lfloor s/r \rfloor +1} \binom{j}{i} p^{(i-2)r}\right)},$$
and so the result now follows from Corollary \ref{linear-zero}. On the other hand, if $l \geq r$ we
have 
$$\langle \zeta^{b} \, | \, \chi \rangle
= \zeta^{b} \sum_{j=0}^{p^{s}-1} \zeta^{p^{2r}\left((up^{l-r}j+b
\sum_{i=2}^{\lfloor s/r \rfloor +1} \binom{j}{i} p^{(i-2)r}\right)}
= \zeta^{b} p^{r} \sum_{j=0}^{p^{s-r}-1} \zeta^{p^{2r}\left((up^{l-r}j+b
\sum_{i=2}^{\lfloor s/r \rfloor +1} \binom{j}{i} p^{(i-2)r}\right)},$$
and so the result now follows by applying Corollary \ref{quadratic-non-zero}
with $k=s-r$, $e=r$ and rescaling appropriately.
\end{proof}

\begin{theorem}\label{cyclotomic-extension}
Let $p$ be an odd prime and let $r \geq 1$. If $p=3$, assume that $r \geq \min\{2,s\}$.
For $k \geq 1$, let
$$R_{k} =
\{ (1+c_{1}p+c_{2}p^{2} + \ldots + c_{\min(r,k)-1}p^{\min(r,k)-1}) \, : 0 \leq c_{i} \leq p-1 \}$$
and let $L_{r+k}$ be the unique subfield of $\Q^{(p^{r+k})}$ with $[\Q^{(p^{r+k})}:L_{r+k}]=p-1$. 
For $s \geq 1$, let
$$ \alpha_{s} = 1 + \sum_{k=1}^{s} \sum_{t \in R_{k}} \tr_{\Q^{(p^{r+k})}/L_{r+k}}(\zeta_{p^{r+k}}^{t})
\quad \textrm{and} \quad G_{s} = \Gal(\Q^{(p^{r+s})}/\Q^{(p^{r})}).
$$
Then for all $s \geq 0$, we have 
$\mathcal{A}_{\Q^{(p^{r+s})}/\Q^{(p^{r})}} = \mathcal{M}_{\Q^{(p^{r})}[G_{s}]}$ 
and $\Z[\zeta_{p^{r+s}}] = \mathcal{A}_{\Q^{(p^{r+s})}/\Q^{(p^{r})}} \cdot \alpha_{s}$. 
In other words, $\Q^{(p^{r+s})}/\Q^{(p^{r})}$ is Leopoldt for all $s \geq 0$.
\end{theorem}

\begin{proof}
That $\mathcal{A}_{\Q^{(p^{r+s})}/\Q^{(p^{r})}} = \mathcal{M}_{\Q^{(p^{r})}[G_{s}]}$ is a 
straightforward application of Theorem \ref{bl-max} (i.e. the corollary to the main theorem in \cite{byott-lettl}).

We prove the main result by induction on $s$. The case $s=0$ is trivial. Now fix $s \geq 1$ and 
assume by induction that $\Z[\zeta_{p^{r+s-1}}] = \mathcal{A}_{\Q^{(p^{r+s-1})}/\Q^{(p^{r})}} \cdot \alpha_{s-1}$.
Henceforth adopt the notation given at the beginning of Section \ref{step}, taking $K=\Q^{(p^{r})}$.
Let $\alpha = \alpha_{s}$ and $\gamma=\alpha_{s-1}$. 
Let $H=\{ h_{1}, \ldots , h_{p-1} \}$ be the unique subgroup of $\Gal(\Q^{(p^{r+s})}/\Q)$ of order $p-1$
and write $h_{i}=\sigma_{i+pf_{i}}$ for some $d_{i} \in \Z$. 
Let $S = \{ (i+pf_{i})x \, | \, x \in R_{s} \}$. Then observe that
$ \alpha = \gamma + \sum_{b \in S} \zeta^{b}.$

We wish to show that the hypotheses of Proposition \ref{check-generator} hold.
Hypothesis (a) is the induction hypothesis and (b) holds by Lemma \ref{pick}.
Thus it remains to check (c).

Let $\chi \in B$. Then $\chi$ is of order $p^{s}$ and there exists a unique $a \in V$ such that
$0 < a < p^{s}$ and $\chi(\sigma_{1+p^{r}})=\zeta_{p^{s}}^{a}=\zeta^{ap^{r}}$.
We can write $a$ uniquely as
$$a=(a_{0}+a_{1}p+a_{2}p^{2}+\ldots+a_{s-1}p^{s-1}), \quad 0 \leq a_{i} \leq p-1, \quad a_{0} \neq 0.$$
Furthermore, there exists a unique element $c \in \N$ of the form
$$ c= (1+c_{1}p+c_{2}p^{2} + \cdots + c_{s-1}p^{s-1}) \quad 0 \leq c_{i} \leq p-1$$
such that for $b:=(a_{0} +pf_{a_{0}})c$, we have
$a \equiv b \, \mod \, p^{s}$, i.e. $v_{p}(b-a) \geq s$.

Suppose that $s \leq r$. Then
$$R_{s} = \{ (1+d_{1}p+d_{2}p^{2} + \cdots + d_{s-1}p^{s-1}) \, : 0 \leq d_{i} \leq p-1 \}$$ 
and so $b$ is the unique element of $S$ such that $v_{p}(b-a) \geq s$. Hence by Lemma \ref{pick}, 
$b$ is the unique element of $S$ such that $\langle \zeta^{b} \, | \, \chi \rangle \neq 0$, 
and so hypothesis (c) of Proposition \ref{check-generator} is satisfied in this case.

Suppose that $s > r$. Then
$$R_{s} = \{ (1+d_{1}p+d_{2}p^{2} + \cdots + d_{r-1}p^{r-1}) \, : 0 \leq d_{i} \leq p-1 \}.$$
There exists a unique $d \in R_{s}$ such that $v_{p}(d-c)\geq r$, namely
$$ d = (1+c_{1}p+c_{2}p^{2} + \cdots + c_{r-1}p^{r-1}).$$
Let $e=(a_{0} +pf_{a_{0}})d \in S$. Since $v_{p}(a_{0} +pf_{a_{0}})=0$, we have
$$v_{p}(e-b) = v_{p}((a_{0}+pf_{a_{0}})(d-c)) = v_{p}(d-c) \geq r.$$
However, $a \equiv b \, \mod \, p^{s}$ and so $v_{p}(e-a) \geq r$.
By construction, $e$ is the unique element of $S$ with this property and therefore
the property that 
$\langle \zeta^{e} \, | \, \chi \rangle \neq 0$ (use Lemma \ref{pick}). Hence hypothesis (c) of 
Proposition \ref{check-generator} is satisfied in this case and so the proof is now complete.
\end{proof}

\begin{corollary}\label{p-power-shifted}
Let $r,s \geq 1$ and $p$ be an odd prime. If $p=3$, assume that $r \geq \min\{2,s\}$.
Let $L/K$ be an extension with $[L:K]=p^{s}$ and $L \subseteq \Q^{(p^{r+s})}$.
Then $\mathcal{A}_{L/K} = \mathcal{M}_{K[G]}$ and $L/K$ is Leopoldt. In fact,
if we suppose that $L$ is of conductor $p^{r+s}$ (which we can do without loss of generality)
then $\mathcal{O}_{L} = \mathcal{A}_{L/K} \cdot \alpha_{s}$ for $\alpha_{s} \in \mathcal{O}_{L}$
as given in Theorem \ref{cyclotomic-extension}.
\end{corollary}

\begin{proof}
That $\mathcal{A}_{L/K} = \mathcal{M}_{K[G]}$ is a straightforward application of Theorem \ref{bl-max} 
(i.e. the corollary to the main theorem in \cite{byott-lettl}).

Assume without loss of generality that $L$ is of conductor $p^{r+s}$.
Let $L'$ and $K'$ be the unique subfields such that $[\Q^{(p^{r+s})}:L']=[\Q^{(p^{r})}:K']=p-1$.
The situation is illustrated by the field diagram below.
$$
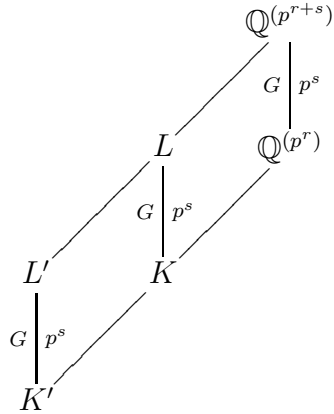
\xymatrix@1@!0@=48pt {
& & \Q^{(p^{r+s})} \ar@{-}[dl]  \ar@{-}[d]_{G}^{p^{s}}  \\
& L  \ar@{-}[dl]  \ar@{-}[d]_{G}^{p^{s}} & \Q^{(p^{r})}  \ar@{-}[dl] \\
 L'  \ar@{-}[d]_{G}^{p^{s}}  &  K  \ar@{-}[dl] & \\
K' & & \\
}
$$
Since their degrees over $K$ are relatively prime, $\Q^{(p^{r})}$ and $L$ are linearly disjoint over $K$. Furthermore, by definition of $\alpha_{s}$ in 
Theorem \ref{cyclotomic-extension}, it is easy to see that in fact 
$\alpha_{s} \in \mathcal{O}_{L'} \subseteq \mathcal{O}_{L}$. Hence the result now follows from
Proposition \ref{shift-down} with $E=\Q^{(p^{r})}$ and $F=\Q^{(p^{r+s})}$.
\end{proof}

\begin{corollary}\label{cyc-Zp}
Let $p$ be an odd prime and $K$ be a subfield of $\Q^{(p)}$. Let 
$$ K=K_{0} \subset K_{1} \subset K_{2} \subset \cdots \subset K_{\infty}^{\mathrm{cyc}} 
$$
be the cyclotomic $\Z_{p}$-extension of $K$. Then for any $n \geq m$, 
the ring of integers $\mathcal{O}_{K_{n}}$ is free as a module over
$\mathcal{A}_{K_{n}/K_{m}}$.
\end{corollary}

\begin{proof}
This is essentially just a rephrasing of Corollary \ref{p-power-shifted}. Recall Remark \ref{exceptional}.
\end{proof}

We can now deduce a theorem that is essentially the same as a result of
Kersten and Michali\v cek (\cite{pNIB}, Theorem 2.1). Also see \cite{greither} and \cite{ichimura}.

\begin{theorem}
Let $p$ be an odd prime and $K$ be any number field. 
Let $$ K=K_{0} \subset K_{1} \subset K_{2} \subset \cdots \subset K_{\infty}^{\mathrm{cyc}} $$
be the cyclotomic $\Z_{p}$-extension of $K$. For every $n \geq m$, the ring of $p$-integers 
$\mathcal{O}_{K_{n}}[p^{-1}]$ is free as a module over 
$\mathcal{O}_{K_{m}}[p^{-1}][\Gal(K_{n}/K_{m})]$.
\end{theorem}

\begin{proof}
The main point is that $p$ is inverted and since only primes above $p$ are ramified in
$K_{n}/K_{m}$, arithmetic disjointness outside of $p$ is equivalent to linear disjointness. Hence the desired
result can be obtained by base change from the case of Corollary \ref{cyc-Zp}.
\end{proof}

\section{Proof of Theorem \ref{maintheorem2}}\label{2-p-power}

We now prove an explicit version of Theorem \ref{maintheorem2}. We adopt the
notation given at the beginning of Section \ref{step}, taking $K=\Q^{(p^{r})+}$.
We abbreviate $\langle - \, | \,  \chi \rangle_{L/K}$ to $\langle - \, | \,  \chi \rangle$. 

\begin{lemma}\label{2pick}
Let $\chi \in B$ and define $a \in \N$ by $\chi (\sigma_{1+p^{r}}) = \zeta_{p^{s}}^{a}=\zeta^{ap^{r}}$.
If $p=3$, assume that $r \geq \min\{2,s\}$.
Let $b \in V$. Then 
$$
\Norm_{\Q^{(p^{r+s})}/\Q}( ( \langle \zeta^{b} \, | \, \chi \rangle) ) = \left\{
\begin{array}{ll}
p^{s(p-1)p^{r+s-1}} & \textrm{if } s \leq r \textrm{ and } \max\{ v_p(b \pm a) \} \geq s \, \\
p^{\frac{1}{2}(s+r)(p-1)p^{r+s-1}} & \textrm{if } s > r \textrm{ and } \max\{ v_p(b \pm a) \} \geq r, \\
0 & \textrm{otherwise.}
\end{array} \right.
$$
\end{lemma}

\begin{proof}
Observe that for any $k \geq 1$, we cannot have both $v_{p}(b-a) \geq k$ and $v_{p}(b+a) \geq k$,
otherwise $v_{p}(b)=v_{p}(2b)=v_{p}((b-a)+(b+a)) \geq k$, contradicting the hypothesis on $b$.
The result now follows immediately from Lemmas \ref{2p-split} and \ref{pick}.
\end{proof}

\begin{theorem}\label{max-real-extension}
Let $p$ be an odd prime and let $r \geq 1$.
If $p=3$, assume that $r \geq \min\{2,s\}$.
For $k \geq 1$, let
$$S_{k} =
\{ (c_{0}+c_{1}p+c_{2}p^{2} + \ldots + c_{\min(r,k)-1}p^{\min(r,k)-1}) \, : 
\, 0 \leq c_{i} \leq (p-1), \,  1 \leq c_{0} \leq \textstyle{\frac{1}{2}}(p-1) \}.$$
For $s \geq 1$, let
$$ \alpha_{s} = \zeta_{p^{r}} + \sum_{k=1}^{s} \sum_{t \in S_{k}} \zeta_{p^{r+k}}^{t}
\quad \textrm{and} \quad G_{s} = \Gal(\Q^{(p^{r+s})}/\Q^{(p^{r})+}).
$$
Then for all $s \geq 0$, we have 
$\Z[\zeta_{p^{r+s}}] = \mathcal{A}_{\Q^{(p^{r+s})}/\Q^{(p^{r})+}} \cdot \alpha_{s}$. 
In other words, $\Q^{(p^{r+s})}/\Q^{(p^{r})+}$ is Leopoldt for all $s \geq 0$.
\end{theorem}

\begin{proof}
We proceed by induction on $s$. In the case $s=0$, it is trivial to check that $\alpha_{0} = \zeta_{p^{r}}$ 
generates a normal integral basis for $\Q^{(p^{r})}/\Q^{(p^{r})+}.$ Now fix $s \geq 1$ and 
assume by induction that $\Z[\zeta_{p^{r+s-1}}] = \mathcal{A}_{\Q^{(p^{r+s-1})}/\Q^{(p^{r})+}} \cdot \alpha_{s-1}$.

Henceforth adopt the notation given at the beginning of Section \ref{step}, taking $K=\Q^{(p^{r})+}$.
Let $\alpha = \alpha_{s}$, $\gamma=\alpha_{s-1}$ and $S=S_{s}$.
Then observe that $ \alpha = \gamma + \sum_{b \in S} \zeta^{b}.$
We wish to show that the hypotheses of Proposition \ref{check-generator} hold.
Hypothesis (a) is the induction hypothesis and (b) holds by Lemma \ref{2pick}.
Thus it remains to check (c).

Let $\chi \in B$. Then $\chi$ is of order $2p^{s}$ and there exists a unique $a \in V$ such that
$0 < a < p^{s}$ and $\chi(\sigma_{1+p^{r}})=\zeta_{p^{s}}^{a}=\zeta^{ap^{r}}$.
We can write $a$ uniquely as
$$a=(a_{0}+a_{1}p+a_{2}p^{2}+\ldots+a_{s-1}p^{s-1}), \quad 0 \leq a_{i} \leq p-1, \quad a_{0} \neq 0.$$

Suppose that $s \leq r$. Then
$$S = \{ (c_{0}+c_{1}p+c_{2}p^{2} + \ldots + c_{s-1}p^{s-1}) \, : 
\, 0 \leq c_{i} \leq (p-1), \,  1 \leq c_{0} \leq \textstyle{\frac{1}{2}}(p-1) \}.$$
Observe that $(S  \biguplus (-S)) \, \mod \, p^{s} \equiv (\Z/p^{s}\Z)^{\times}$ and so there exists a unique
element $b \in S$ such that either $v_{p}(b-a) \geq s$ or $v_{p}(b+a) \geq s$.
Hence by Lemma \ref{2pick}, $b$ is the unique element of $S$ such that 
$\langle \zeta^{b} \, | \, \chi \rangle \neq 0$, 
and so hypothesis (c) of Proposition \ref{check-generator} is satisfied in this case.

Suppose that $s > r$. 
Then $$S = \{ (c_{0}+c_{1}p+c_{2}p^{2} + \ldots + c_{r-1}p^{r-1}) \, : 
\, 0 \leq c_{i} \leq (p-1), \,  1 \leq c_{0} \leq \textstyle{\frac{1}{2}}(p-1) \}.$$
Observe that $(S \biguplus (-S)) \, \mod \, p^{r} \equiv (\Z/p^{r}\Z)^{\times}$ and so
there exists a unique $b \in S$ such that
$$ \pm b \equiv (a_{0} + a_{1}p + \ldots + a_{r-1}p^{r-1}) \, \mod \, p^{r}.$$
This is the unique element of $S$ such that either 
$v_{p}(b-a) \geq r$ or $v_{p}(b+a) \geq r$ and so the unique
element of $S$ such that $\langle \zeta^{b} \, | \, \chi \rangle \neq 0$ (use Lemma \ref{2pick}).
Hence hypothesis (c) of Proposition \ref{check-generator} is 
satisfied in this case and so the proof is now complete.
\end{proof}

\begin{remark}
The idempotent
$ e = \frac{1}{2p^{s}}\tr_{\Q^{(p^{r+s})}/\Q^{(p^{r})+}}$
is in $\mathcal{M}_{\Q^{(p^{r})+}[G_{s}]}$, but is not in 
$\mathcal{A}_{\Q^{(p^{r+s})}/\Q^{(p^{r})+}}$ 
because $e(\zeta_{p^{r}})=\frac{1}{2}(\zeta_{p^{r}}+\zeta_{p^{r}}^{-1}) \notin \Z[\zeta_{p^{r+s}}].$
Therefore $\mathcal{A}_{\Q^{(p^{r+s})}/\Q^{(p^{r})+}} \subsetneqq \mathcal{M}_{\Q^{(p^{r})+}[G_{s}]}$. 
\end{remark}

\section{Acknowledgments}

The author is grateful to David Burns, Steven Chase, Griff Elder, Eknath Ghate, 
Gerhard Michler, Shankar Sen and Anupam Srivastav for useful conversations; 
to Nigel Byott, Spencer Hamblen, \'Alvaro Lozano-Robledo and Jason Martin for looking at various 
drafts of this paper; and to Humio Ichimura for providing off-prints of several of his papers. 
The author also wishes to express his gratitude to Cornelius Greither, Ravi Ramakrishna and the referee for numerous helpful comments and suggestions.

\nocite{jaco}
\nocite{frohlich-book}
\bibliography{LeoAccBib}{}
\bibliographystyle{amsalpha}

\end{document}